\documentclass[12pt]{article}
\usepackage{amsfonts}

\usepackage{latexsym}
\usepackage{amssymb}
\usepackage{epsfig}
\usepackage{amsmath}
\usepackage{amsthm}
\usepackage[mathscr]{eucal}
\usepackage{subfigure}
\usepackage{graphicx}
\usepackage{hyperref}

\newtheorem{Lemma}{Lemma}

\newtheorem{Theorem}[Lemma]{Theorem}

\newtheorem{Corollary}[Lemma]{Corollary}

\renewcommand{\qed}{\hfill{\ \ \rule{2mm}{2mm}} \vspace{0.2in}}

\newcommand{\ind}{1\hspace{-2.3mm}{1}}

\begin{document}

\title{Long Paths and Hamiltonian paths in Inhomogenous Random Graphs}
\author{ \textbf{Ghurumuruhan Ganesan}
\thanks{E-Mail: \texttt{gganesan82@gmail.com} } \\
\ \\
New York University, Abu Dhabi }
\date{}
\maketitle

\begin{abstract}
In this paper, we study long paths and Hamiltonian paths in inhomogenous random graphs. In the first part of the paper,
we consider an inhomogenous Erd\H{o}s-R\'enyi random graph~\(G_E\) with average edge density~\(p_n.\)
We prove that if~\(np_n^2 \longrightarrow \infty\) as~\(n \rightarrow \infty,\) then the longest path contains
at least~\(n-ne^{-\delta_1 np_n^2}\) nodes with high probability (i.e., with probability converging to
one as~\(n \rightarrow \infty\)), for some constant~\(\delta_1> 0 .\) In particular, if~\(np_n^2 = M\log{n}\)
for some constant~\(M > 0\) large, then~\(G_E\) is Hamiltonian with high probability; i.e., the longest path
contains all the nodes of~\(G_E.\)

In the second part of the paper, we consider a random geometric graph~\(G_R\) consisting of~\(n\) nodes,
each independently distributed according to a (not necessarily uniform) density~\(f.\)
If~\(r_n\) is the connectivity radius and~\(nr_n^2 \longrightarrow \infty,\) then with high probability,
the longest cycle contains at least~\(n-ne^{-\delta_2 nr_n^2}\) nodes for some constant~\(\delta_2 > 0.\)
As a consequence of our proof, we obtain that if~\(nr_n^2 = \log{n} + 7\log{\log{n}} + \omega_n\)
and~\(\omega_n \longrightarrow \infty\) as~\(n \rightarrow \infty,\) then with high probability~\(G_R\) contains
a Hamiltonian cycle.


\vspace{0.1in} \noindent \textbf{Key words:} Inhomogenous random graphs, random geometric graphs, long paths, Hamiltonian paths.

\vspace{0.1in} \noindent \textbf{AMS 2000 Subject Classification:} Primary:
60J10, 60K35; Secondary: 60C05, 62E10, 90B15, 91D30.
\end{abstract}

\bigskip

\setcounter{equation}{0}
\renewcommand\theequation{\thesection.\arabic{equation}}
\section{Introduction}

\subsection{Erd\H{o}s-R\'enyi (ER) Random Graphs}
Let~\(K_n\) be the labelled complete graph on~\(n \geq 3\) vertices with vertex set\\\(\{1,2\ldots,n\}\) and let~\(e(i,j)\) denote the edge between vertices~\(i\) and~\(j.\) Let~\(X(i,j)\) be a Bernoulli random variable defined on the probability space\\\((\{0,1\}, \mathbb{B}(\{0,1\}), \mathbb{P}_{i,j})\) with~\[\mathbb{P}_{i,j}(X(i,j) = 1) = p(i,j) = 1 -\mathbb{P}_{i,j}(X(i,j) = 0).\] Here~\(\mathbb{B}(\{0,1\})\) is the set of all subsets of~\(\{0,1\}.\) We say that edge~\(e(i,j)\) is \emph{open} if~\(X(i,j) = 1\) and closed otherwise. The random variables~\(\{X(i,j)\}\) are independent and the resulting random graph~\(G\) is an inhomogenous Erd\H{o}s-R\'enyi (ER) random graph, defined on the probability space~\((\Omega,{\cal F}, \mathbb{P}).\) Here~\(\Omega = \{0,1\}^{ {n \choose 2}},\) the sigma algebra~\({\cal F}\) is the set of subsets of~\(\Omega\) and~\(\mathbb{P} = \prod_{i,j} \mathbb{P}_{i,j}.\)

We assume that there is a sequence~\(p_n \in (0,1), n \geq 3\) and constants~\(0 < \beta_1 \leq \beta_2 < \beta_3 \leq 1\) so that
\begin{equation}\label{beta_3def}
\inf_{1 \leq i \leq n} \frac{1}{n-1}\sum_{j\neq i} p(i,j) \geq \beta_3 p_n
\end{equation}
and
\begin{equation}\label{beta_12def}
\inf_{1 \leq i \leq n}\inf_{S} \frac{1}{\#S}\sum_{j \in S} p(i,j) \geq \beta_1 p_n
\end{equation}
for all~\(n\) large. For a fixed~\(1 \leq i \leq n,\) the infimum above is taken over all sets~\(S\) such that~\(\#S \geq \beta_2 n p_n\) and~\(i \notin S.\) The first condition~(\ref{beta_3def}) implies that the average number of neighbours per vertex is at least~\(\beta_3 (n-1)p_n\) and the second condition implies that the average edge density taken over sets of cardinality at least~\(\beta_2 n p_n,\) is at least~\(\beta_1 p_n.\) All constants mentioned are independent of~\(n.\)

We have the following result.
\begin{Theorem}\label{long_er_inhom} Suppose
\begin{equation}\label{pn_cond}
p_n \longrightarrow 0 \text{ and } \frac{\log{n}}{np_n} \longrightarrow 0
\end{equation}
as~\(n \rightarrow \infty.\) If~\(L_n\) denotes the length of the longest path in the random graph~\(G = G(n,p_n),\) then
\begin{equation}\label{exp_ln}
\mathbb{E}L_n \geq n-2ne^{-\beta_1 \beta_2 np_n^2}
\end{equation}
for all~\(n\) large.
Also
for any~\(0 < \delta < 1,\) we have
\begin{equation}\label{eq_lg1}
\mathbb{P}\left(L_n \geq  n - 2n\exp\left(-\beta_1\beta_2(1-\delta)np_n^2\right)\right) \geq 1 - \exp\left(-\beta_1\beta_2\delta np_n^2\right)
\end{equation}
for all~\(n\) large. Suppose that~\(np_n^2 = M\log{n}\) for some~\(M > (\beta_1\beta_2)^{-1}.\) Setting~\(M_1 = M\beta_1\beta_2 > 1,\) we have
\begin{equation}\label{expl_ln2}
\mathbb{E}L_n \geq n - \frac{2}{n^{M_1-1}}
\end{equation}
and
\begin{equation}\label{eq_lg2}
\mathbb{P}(L_n = n) \geq 1 - \frac{2}{n^{M_1-1}}
\end{equation}
for all~\(n\) large.
\end{Theorem}
The final result implies that the random graph~\(G\) contains a Hamiltonian path with high probability. For homogenous random graphs, the standard methods to study long paths and Hamiltonian paths usually include a combination of edge sprinkling, Markov chain analysis and path rotation (see Bollobas~(2001), Chapter~\(8\) and references therein). For inhomogenous graphs as described above, the above methodology is not directly applicable since the individual edge probabilities could be arbitrarily low. We use a simple subgraph analysis technique to study the long paths in inhomogenous random graphs~(see Section~\ref{pf_long_path_er2}).

For homogenous ER graphs, we in fact have the following Corollary.
\begin{Corollary}\label{long_er_hom}
Suppose~(\ref{pn_cond}) holds. For any~\(0 < \delta < \frac{1}{2}\) we have
\begin{equation}\label{exp_ln_hom}
\mathbb{E}L_n \geq n-2ne^{-(1-\delta)np_n^2}
\end{equation}
and
\begin{equation}\label{eq_lg1_hom}
\mathbb{P}\left(L_n \geq  n - 2ne^{-(1-2\delta)np_n^2}\right) \geq 1 - e^{-\delta np_n^2}
\end{equation}
for all~\(n\) large.
Suppose that~\(np_n^2 = M\log{n}\) for some~\(M > 1.\) For any~\(1 < M_1 < M,\) we have
\begin{equation}\label{exp_ln2_hom}
\mathbb{E}L_n \geq 1 - \frac{2}{n^{M_1-1}}
\end{equation}
and
\begin{equation}\label{eq_lg2_hom}
\mathbb{P}(L_n = n) \geq 1 - \frac{2}{n^{M_1-1}}
\end{equation}
for all~\(n\) large.
\end{Corollary}


\subsection{Random Geometric Graphs}
Consider $n$ vertices \(X_1,X_2,\ldots,X_n,\) independently distributed in the unit square $S = \left[-\frac{1}{2},\frac{1}{2}\right]^2$ each according to a certain density \(f\) satisfying
\begin{equation}\label{f_eq}
0 < \inf_{x \in S} f(x) \leq \sup_{x \in S} f(x) < \infty.
\end{equation}
We define the overall process on the probability space \((\Omega_X, {\cal F}_X, \mathbb{P}).\) Connect two vertices \(X_i\) and \(X_j\) by a edge \(e\) if the Euclidean distance $d(X_i,X_j)$ between them is less than $r_n.$ The resulting graph is denoted as $G = G(n,r_n,f)$ and called a random geometric graph~(RGG). Let \(C_G\) denote the component of \(G\) containing the largest number of nodes. In Ganesan (2013), we have proved that if \(nr_n^2 \longrightarrow \infty,\) then with high probability (i.e. with probability tending to one as \(n \rightarrow \infty\)), the largest component \(C_G\) contains at least \(n - ne^{-\beta nr_n^2}\) nodes, for some constant \(\beta > 0.\) 

In this paper we study the number of edges in the longest cycle in the random graph~\(G.\) Let \(LC_n\) denote the length of the longest cycle in~\(G_n = G(n,r_n,f).\) We have the following result.
\begin{Theorem}\label{long_rgg}
Suppose~
\begin{equation}\label{assumption}
nr_n^2 \longrightarrow \infty
\end{equation}
as~\(n \rightarrow \infty.\) There are constants~\(\delta_1,\delta_2> 0\) such that
\begin{equation}\label{eq_lg1_rgg}
\mathbb{P}\left(LC_n \geq  n - ne^{-\delta_1 nr_n^2}\right) \geq 1 - e^{-\delta_2 nr_n^2}
\end{equation}
for all~\(n\) large. Suppose that~\(f\) is uniform and~
\begin{equation}\label{rn_hom}
nr_n^2 = \log{n} + 7\log{\log{n}} + \omega_n
\end{equation}
where~\(\omega_n \rightarrow \infty\) as~\(n \rightarrow \infty.\) We have
\begin{equation}\label{eq_lg2_rgg}
\mathbb{P}(LC_n = n) \geq 1 - Ce^{-\omega_n}
\end{equation}
for all~\(n\) large and for some constant~\(C > 0.\)
\end{Theorem}
The first result~(\ref{eq_lg1_rgg}) obtains estimates on the length of long cycles for the subconnective case where~\(nr_n^2 \longrightarrow \infty.\) This extends previous results~(see Diaz et al (2007), Balogh et al (2011) and references therein) which have primarily studied Hamiltonian cycles in RGGs in the connectivity regime. In fact, as a corollary of our proof technique, we also obtain~(\ref{eq_lg2_rgg}) stating that slightly above the connectivity regime, the random geometric graph contains a Hamiltonian cycle with high probability.

The paper is organized as follows. In Section~\ref{pf_long_path_er2}, we prove Theorem~\ref{long_er_inhom} and obtain Corollary~\ref{long_er_hom} as a Corollary. In Section~\ref{pf_long_cycle_rgg}, we prove estimate~(\ref{eq_lg1_rgg}) in Theorem~\ref{long_rgg}. Finally, in Section~\ref{pf_long_cycle_rgg_ham}, we prove estimate~(\ref{eq_lg2_rgg}) in Theorem~\ref{long_rgg}.

\setcounter{equation}{0}
\renewcommand\theequation{\thesection.\arabic{equation}}
\section{Proof of Theorem~\ref{long_er_inhom}}\label{pf_long_path_er2}
Let~\(G\) denote the random graph~\(G(n,p_n)\) and let~\(P(G)\) be the (random) longest path in~\(G.\) If there is more than one choice, we choose one according to some predetermined order, e.g., lexicographic ordering. Both the estimates~(\ref{eq_lg1}) and~(\ref{eq_lg2}) follow from the below estimate. We have
\begin{eqnarray}
\sup_{1 \leq i \leq n} \mathbb{P}\left(i \notin P(G)\right) \leq 2\exp\left(- \beta_1 \beta_2  np_n^2\right)\label{eq_fin_lg1_est_c2}
\end{eqnarray}
for all~\(n\) large. 

To prove~(\ref{eq_fin_lg1_est_c2}), we need a preliminary estimate. Fix~\(1 \leq i \leq n\) and let~\(G_i\) denote the random induced subgraph formed by the vertices~\(\{1,\ldots,n\}\setminus \{i\}.\) Fix~\(\delta > 0\) small and let~\(A_i\) denote the event that every vertex in the graph~\(G_i\) has degree at least
\begin{equation}\label{t_0_def2}
t_0 := \beta_2 np_n
\end{equation}
for all~\(n\) large. Here~\(\beta_2 \leq 1\) is as in~(\ref{beta_12def}). Using Chernoff bounds, we have the following estimate for the event~\(A_i.\)
\begin{equation}\label{a1_est_eq_c2}
\sup_{1 \leq i \leq n} \mathbb{P}(A^c_i)  \leq a_n := n\exp\left(-q(\delta)\beta_2 np_n\right)
\end{equation}
for all~\(n\) large, where~\(q(\delta) > 0\) satisfies
\begin{equation}\label{qalp_def}
e^{-q(\delta)} = \min\left(\frac{e^{\delta}}{(1+\delta)^{1+\delta}},\frac{e^{-\delta}}{(1-\delta)^{1-\delta}}\right).
\end{equation}
\emph{Proof of~(\ref{a1_est_eq_c2})}:  We use the following Chernoff bound. Let~\(\{X_j\}_{1 \leq j \leq m}\) be independent Bernoulli random variables with~\[\mathbb{P}(X_j = 1) = p_j = 1-\mathbb{P}(X_j = 0).\] We have the following estimate. 
Fix~\(\alpha > 0.\) If~\[T_m = \sum_{j=1}^{m} X_j\] and~\(\mu_m = \mathbb{E}T_m,\) then
\begin{equation}\label{conc_est_f}
\mathbb{P}\left(\left|T_m - \mu_m\right| \geq \mu_m \alpha \right) \leq 2\exp\left(-q(\alpha)\mu_m\right)
\end{equation}
for all \(m \geq 1,\) where~\(q(\alpha)\) is as in~(\ref{qalp_def}). The above result follows using Chernoff bounds (for a proof, we refer to the Wikipedia link\\\(https://en.wikipedia.org/wiki/Chernoff\_bound\)). 

Let~\(E(j)\) denote the event that vertex~\(j\) has at least~\(t_0\) neighbours in the random graph~\(G_i\) so that~\(A_i = \bigcap_{j=1, j\neq i}^{n} E(j)\) and
\begin{equation}\label{a1_def2}
\mathbb{P}(A^c_i) \leq \sum_{j \neq i}\mathbb{P}(E^c(j))
\end{equation}
Fixing~\(j \neq i ,\) we estimate each~\(E^c(j)\) separately. Let
\begin{equation}\label{mu_i_def}
\mu_n(j) = \sum_{k=2}^{j-1} p(j,k) + \sum_{k=j+1}^{n} p(j,k)
\end{equation}
be the mean number of neighbours of vertex~\(j\) in the graph~\(G_i.\) Using~(\ref{beta_3def}), we have
\begin{equation}\label{mu_i_est}
\mu_n(j) \geq \beta_3(n-1)p_n-1
\end{equation}
for all~\(n \geq N_1.\) Here~\(N_1 \geq 1\) does not depend on~\(i\) or~\(j.\) We recall that~\(\beta_2 < \beta_3\) (see~(\ref{beta_12def})) and so choosing~\(\delta > 0\) small so that~\(\beta_3(1-\delta) > \beta_2,\) we have
\begin{equation}\label{mu_est2}
\mu_n(j)(1-\delta) \geq (\beta_3(n-1)p_n-1)(1-\delta)  \geq \beta_2 np_n  = t_0
\end{equation}
for all~\(n \geq N_2\) large, where the final equality follows from~(\ref{t_0_def2}). Here~\(N_2 = N_2(\beta_2,\beta_3,\delta)\) does not depend in~\(i\) or~\(j.\)



Using the estimate~(\ref{conc_est_f}) with~\(m= n-2,\mu_m = \mu_n(j)\) and~\(\alpha = \delta,\) we have
\begin{eqnarray}\label{ei_est11}
\mathbb{P}(E^c(j)) \leq \exp\left(-q(\delta)\mu_n(j)\right) \leq \exp\left(-q(\delta)\beta_2 np_n\right).
\end{eqnarray}
for all~\(n \geq N_2.\) Here~\(N_2 \geq 1\) is as in~(\ref{mu_est2}).
Using~(\ref{ei_est11}) in~(\ref{a1_def2}), we have
\begin{equation}\label{a1_est2}
\mathbb{P}(A_i^c) \leq n\exp\left(-q(\delta)\beta_2 np_n\right)
\end{equation}
for all~\(n \geq N_2.\) Since~\(N_2\) does not depend on~\(i,\) this proves~(\ref{a1_est_eq_c2}).~\(\qed\)

We use~(\ref{a1_est_eq_c2}) to prove~(\ref{eq_fin_lg1_est_c2}).\\
\emph{Proof of~(\ref{eq_fin_lg1_est_c2})}: We have
\begin{eqnarray}
\mathbb{P}\left(i \notin P(G)\right) &=& \mathbb{P}\left(\left\{i \notin P(G)\right\} \bigcap A_i\right) +
\mathbb{P}\left(\left\{i \notin P(G)\right\} \bigcap A_i^c\right)\nonumber\\
&\leq& \mathbb{P}\left(\left\{i \notin P(G)\right\} \bigcap A_i\right) + \mathbb{P}(A_i^c) \nonumber\\
&\leq& \mathbb{P}\left(\left\{i \notin P(G)\right\} \bigcap A_i\right) + a_n\label{1pg_1_c2}
\end{eqnarray}
where the sequence~\(a_n\) is as defined in~(\ref{a1_est_eq_c2}).

Suppose now that the event~\(\{i \notin P(G)\} \bigcap A_i\) occurs. Since the vertex~\(i\) does not belong to~\(P(G),\) the longest path in the graph~\(G_i\) is also~\(P(G);\) i.e., the path~\(P(G_i) = P(G).\) We therefore have
\begin{eqnarray}
&&\mathbb{P}\left(\left\{i \notin P(G)\right\} \bigcap A_i\right) \nonumber\\
&&\;\;\;= \mathbb{P}\left(\left\{i \notin P(G)\right\} \bigcap A_i \bigcap \{P(G_i) = P(G)\}\right) \nonumber\\
&&\;\;\;= \sum_{\pi} \mathbb{P}\left(\{P(G) = \pi\} \bigcap \left\{i \notin \pi\right\} \bigcap A_i \bigcap \{P(G_i) = \pi\} \right)\label{1pg_2_c2}
\end{eqnarray}
where the summation is taken over all paths~\(\pi\) formed by the vertices~\(\{1,\ldots,n\} \setminus \{i\}.\)

For a fixed path~\(\pi = (\pi(1),\ldots,\pi(f)),\) we let~\(f = \#\pi\) denote the number of vertices of~\(\pi\) and let~\(\pi(1)\) denote the least endvertex of~\(\pi;\) i.e., \(\pi(1) < \pi(f).\) For a fixed~\(\pi,\) suppose that the event in the brackets in the right hand side of~(\ref{1pg_2_c2}) occurs. Let~\(N(\pi(1))\) denote the set of neighbours of the endvertex~\(\pi(1)\) in the graph~\(G_i.\) We therefore have
\begin{eqnarray}
&&\mathbb{P}\left(\{P(G) = \pi\} \bigcap \left\{i \notin \pi\right\} \bigcap A_i \bigcap \{P(G_i) = \pi\} \right) \nonumber\\
&&\;\;\;= \sum_{S} \mathbb{P}\left(\{P(G) = \pi\} \bigcap \left\{i \notin \pi\right\} \bigcap A_i \bigcap T(\pi,S) \right) \label{1pg_3_c2}
\end{eqnarray}
where the event
\begin{equation}\label{t_pi_s_def}
T(\pi,S) = \{P(G_i) = \pi\} \bigcap \{N(\pi(1)) = S\}
\end{equation}
and the summation is over all subsets of~\(\{1,\ldots,n\} \setminus \{i\}.\) Fix a set~\(S\) and suppose that the event within the brackets of the final term in~(\ref{1pg_3_c2}) occurs. We have the following properties.\\
\((a1)\) All the neighbours of the endvertex~\(\pi(1)\) belong to the path~\(\pi;\) in other words, the set~\[S = \{\pi(j_1),\ldots,\pi(j_t)\}\] for some indices~\(1 \leq j_1 \leq \ldots \leq j_t \leq f.\) \\
\((a2)\) The set~\(S\) contains~\(t \geq t_0\) vertices, where~\(t_0\) is as defined in~(\ref{t_0_def2}).\\
\((a3)\) The vertex~\(i\) is not adjacent to any of the vertices in the set~\[R = \{\pi(j_1-1),\ldots,\pi(j_t-1)\}.\]

\begin{figure}[tbp]
\centering
\includegraphics[width=2.5in, trim= 100 220 80 260, clip=true]{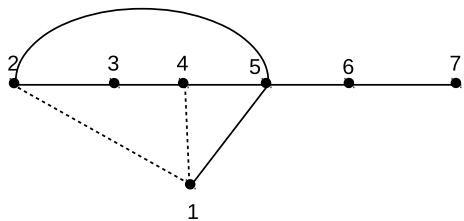}
\caption{The longest path~\(P(G)\) in~\(G\) is~\(\pi = (2,3,4,5,6,7).\) The vertex~\(\pi(1) = 2\) and the neighbour set~\(N(\pi(1)) = S = \{3,5\}.\) Since~\(1 \notin P(G),\) the vertex~\(1\) cannot be adjacent any vertex in the set~\(R = \{2,4\}.\) If for example~\(1\) were adjacent to~\(4,\) then we would have a longer path in~\(G\) formed by~\((1,4,3,2,5,6,7).\)}
\label{fig_long}
\end{figure}

The property~\((a3)\) is illustrated in Figure~\ref{fig_long} where~\(n = 7\) and the longest path in the random graph~\(G\) is given by path \(P(G) = \pi = (2,3,4,5,6,7).\) The vertex~\(1 \notin P(G) = \pi\) and~\(\pi(1) =2\) and~\(N(\pi(1)) = S = \{3,5\}.\) The vertex~\(1\) is not adjacent to any vertex in the set~\(R = \{2,4\}.\) If for example~\(1\) were adjacent to~\(4,\) then~\((1,4,3,2,5,6,7)\) would form a longer path in~\(G.\)\\\\
\emph{Proof of~\((a1)-(a3)\)}: The property~\((a1)\) is true since~\(\pi\) is the longest path in the graph~\(G_i.\) If~\(\pi(1)\) contains a neighbour~\(z \notin \pi\) in the graph~\(G_i,\) then~\((z,\pi(1),\ldots,\pi(f))\) would form a longer path in~\(G_i.\) The property~\((a2)\) is true since the event~\(A_i\) occurs~(see paragraph preceding~(\ref{t_0_def2})) and so every vertex in~\(G_i\) has at least~\(t_0\) neighbours in~\(G_i.\)

To prove~\((a3),\) we use the fact that the event~\(\{P(G) = \pi\}\cap \{i \notin \pi\}\) occurs. So the vertex~\(i\) is not adjacent to any of the vertices in the set~\(R.\)  Because, otherwise, we would obtain a path of longer length in~\(G.\) For example, if~\(i\) was adjacent to~\(\pi(j_1-1),\) then~\[(i,\pi(j_1-1),\pi(j_1-2),\ldots,\pi(1),\pi(j_1),\pi(j_1+1),\ldots,\pi(f))\] would form a path in~\(G\) containing one more edge than~\(\pi;\) i.e., \(f\) edges. This contradicts the fact that the event~\(\{P(G) = \pi\}\) occurs and so every path in~\(G\) has at most~\(f-1\) edges.~\(\qed\)

From property~\((a2)\) above we have that the set~\(R\) contains at least~\(t_0\) vertices and so define~\(V(S)\) to be the event that the vertex~\(i\) is not adjacent to any of the vertices in~\(\{\pi(j_1-1),\ldots,\pi(j_{t_0-1})\}.\)
From property~\((a3),\) we therefore have
\begin{eqnarray}
&&\mathbb{P}\left(\{P(G) = \pi\} \bigcap \left\{i \notin \pi\right\} \bigcap A_i \bigcap T(\pi,S) \right)\nonumber\\
&&\;\;\;\leq \mathbb{P}\left(T(\pi,S)\bigcap V(S)\right) \nonumber\\
&&\;\;\;= \mathbb{P}\left(T(\pi,S)\right)\mathbb{P}(V(S)). \label{1pg_4_c2}
\end{eqnarray}
The equality~(\ref{1pg_4_c2}) is true as follows. We recall that the event~\(T(\pi,S)\) (see~(\ref{t_pi_s_def})) depends only on the state of edges with vertices in the graph~\(G_i\) and from the definition above, the event~\(V(S)\) depends on the state of edges containing~\(i\) as an endvertex. Therefore the events~\(T(\pi,S)\) and~\(V(S)\) are independent.

For a fixed set~\(S\) we have the following estimate for the event~\(V(S).\) Letting~\(\beta_1,\beta_2 > 0\) be as in~(\ref{beta_12def})
we have
\begin{equation}\label{vs_est_c2}
\mathbb{P}(V(S)) \leq v_n := \exp\left(-\beta_1\beta_2np^2_n\right)
\end{equation}
for all~\(n \geq N\) large. Here~\(N = N(\beta_1,\beta_2)\) does not depend on the choice of~\(S.\)\\
\emph{Proof of~(\ref{vs_est_c2})}:
Since the set~\(R\) contains \[\#R \geq t \geq t_0 = \beta_2np_n\]  vertices (see~(\ref{t_0_def2})), we have using~(\ref{beta_12def}) that
\begin{equation}\nonumber
\mathbb{P}(V(S)) \leq \prod_{j \in R}(1-p(i,j)) \leq \exp\left(-\sum_{j\in R} p(i,j)\right) \leq \exp\left(-\beta_1 \beta_2 np^2_n\right)
\end{equation}
for all~\(n\) large.~\(\qed\)

Substituting~(\ref{vs_est_c2}) into~(\ref{1pg_4_c2}) and using~(\ref{1pg_3_c2}) we get
\begin{eqnarray}
&&\mathbb{P}\left(\{P(G) = \pi\} \bigcap\left\{i \notin \pi\right\} \bigcap A_i \bigcap \{P(G_i) = \pi\}\right) \nonumber\\
&&\;\;\;\leq \sum_{S} \mathbb{P}\left(T(\pi,S)\right) v_n \label{1pg_33}
\end{eqnarray}
and substituting the above into~(\ref{1pg_2_c2}) gives
\begin{eqnarray}
\mathbb{P}\left(\left\{i \notin P(G)\right\} \bigcap A_i\right) \leq \left(\sum_{\pi}\sum_S \mathbb{P}\left(T(\pi,S)\right) \right) v_n.\nonumber
\end{eqnarray}
The events~\(T(\pi,S)\) are disjoint for distinct pairs~\((\pi,S)\) and so we have
\begin{eqnarray}
\mathbb{P}\left(\left\{i \notin P(G)\right\} \bigcap A_i\right) \leq v_n.\label{1pg_11}
\end{eqnarray}
Substituting the above into~(\ref{1pg_1_c2}), we have
\begin{eqnarray}
\mathbb{P}\left(i \notin P(G)\right) &\leq& v_n + a_n \label{van_est}\\
&\leq& \exp\left(-\beta_1\beta_2 np^2_n\right) + n\exp\left(-q(\delta)\beta_2 np_n\right)\nonumber\\
&\leq& 2\exp\left(-\beta_1 \beta_2 np_n^2\right) \label{van_est2}
\end{eqnarray}
for all~\(n\) large. To see that the final estimate is true, it is enough to see that
\[\beta_1 \beta_2 np_n^2 < q(\delta) \beta_2 np_n - \log{n}\] for all~\(n\) large. Equivalently, it is enough to see that \[\beta_1\beta_2 p_n < q(\delta)\beta_2 - \frac{\log{n}}{np_n}\] for all~\(n\) large, which is true since~\(p_n \rightarrow 0\) and~\(\frac{\log{n}}{np_n} \longrightarrow 0\) as~\(n \rightarrow \infty\)~(see~(\ref{pn_cond})). This proves~(\ref{van_est2}).

Using~(\ref{eq_fin_lg1_est_c2}), we obtain~(\ref{eq_lg1}) and~(\ref{eq_lg2}) as follows. Let~\[X_O = X_O(G) := \sum_{j=1}^{n} \ind(j \notin P(G))\] denote the set of vertices not belonging to the longest path~\(P(G)\) in~\(G.\) We have from the estimate~(\ref{eq_fin_lg1_est_c2}) that
\begin{equation}\label{exo_est}
\mathbb{E}X_O \leq 2ne^{-\beta_1 \beta_2np_n^2}.
\end{equation}
This proves~(\ref{exp_ln}) and using Markov inequality, we have
\[\mathbb{P}\left(X_O \geq 2n\exp\left(-\beta_1\beta_2(1-\delta)np_n^2\right)\right) \leq \exp\left(-\beta_1\beta_2\delta np_n^2\right)\] for any~\(0 < \delta < 1.\) This proves~(\ref{eq_lg1}).

To prove~(\ref{eq_lg2}), we assume that~\(np_n^2 = M\log{n}\) for some constant~\(M > (\beta_1\beta_2)^{-1}.\) We then obtain from~(\ref{exo_est}) that~\[\mathbb{E}X_O \leq 2ne^{-\beta_1\beta_2 np_n^2} \leq \frac{2}{n^{M_1-1}}\] where~\(M_1 = M\beta_1\beta_2 > 1.\) This proves~(\ref{expl_ln2}) and again using Markov inequality, we have \[\mathbb{P}(X_O \geq 1) \leq \mathbb{E}X_O \leq \frac{2}{n^{M_1-1}}.\] This proves~(\ref{eq_lg2}).~\(\qed\)

\emph{Proof of Corollary~\ref{long_er_hom}}: Here~(\ref{beta_3def}) and~(\ref{beta_12def}) are satisfied with~\(\beta_3 = \beta_1 = 1.\) And so the estimate for the sequences~\(a_n\) and~\(v_n\) in~(\ref{a1_est_eq_c2}) and~(\ref{vs_est_c2}) hold with~\(\beta_1 = \beta_2 = 1.\)~\(\qed\)


\section{Proof of~(\ref{eq_lg1_rgg}) in Theorem~\ref{long_rgg}}\label{pf_long_cycle_rgg}
For integer \(n \geq 1,\) let
\begin{equation}\label{k_n_def}
K_n := \left\lceil\frac{\log{n}}{nr_n^2}\right\rceil
\end{equation}
where \(\lceil x \rceil\) refers to the smallest integer strictly larger than \(x.\) We need the following estimate for future use. For all \(n\) large, we have
\begin{equation}\label{k_n_est}
\frac{1}{\sqrt{n}} \leq r_n \leq K_n r_n \leq K_n^2r_n \leq \max\left(4\frac{(\log{n})^{2.5}}{\sqrt{n}}, r_n\right) \longrightarrow 0
\end{equation}
as \(n \rightarrow \infty\) where the final convergence follows from  (\ref{assumption}). We use (\ref{assumption}) to get that \(nr_n^2 \geq 1\) for all \(n\) large. This proves the first inequality. The second and the third inequalities are obtained using \(K_n \geq 1\) for all \(n \geq 1.\) We obtain the final inequality as follows. If \( nr_n^2 \leq \log{n},\) then \(r_n \leq \sqrt{\frac{\log{n}}{n}}\) and we use \(nr_n^2 \geq 1\) to get that \( K_n  = \lceil \frac{\log{n}}{nr_n^2} \rceil \leq 2\log{n}\) for all \(n\) large. This implies that \(K_n^2 r_n \leq 4\frac{(\log{n})^{2.5}}{\sqrt{n}}.\) If \(nr_n^2 \geq \log{n}, \) we have \(K_n = 1\) and so \(K_n^2 r_n = r_n.\)

\subsection*{Construction of the backbone}
Tile the unit square~\(S\) into disjoint squares~\(\{S_j\}\)  each of size~\(\epsilon_1 r_n \times \epsilon_1 r_n.\) Here~\(\epsilon_1 = \epsilon_1(n) \in \left(\frac{1}{4},\frac{1}{5}\right)\) so that~\(\frac{1}{\epsilon_1 r_n}\) is an integer. We choose~\(\epsilon_1\) as above so that the following condition holds: If~\(S_{i_1}\) and~\(S_{i_2}\) are two squares which share a corner, then every node in~\(S_{i_1}\) is connected to every node in~\(S_{i_2}\) by an edge.

Divide the unit square~\(S\) into a set of horizontal rectangles \({\cal R}_H\) each of size $1 \times MK_n\epsilon_1 r_n$ and also vertically into a set of rectangles \({\cal R}_V,\) each of size $MK_n \epsilon_1 r_n \times 1.$ If \((MK_n\epsilon_1 r_n)^{-1}\) is an integer, we obtain a perfect tiling as in Figure~3(a) of Ganesan (2013). Otherwise we start the tiling from the bottom until we reach close to the top and add another \(1 \times MK_n\epsilon_1 r_n\) rectangle sharing the top edge with~\(S.\) Thus the two top most rectangles in the tiling overlap as in Figure~3(b) of Ganesan~(2013).

\begin{figure}[tbp]
\centering
\subfigure{
   \includegraphics[width = 2.5in] {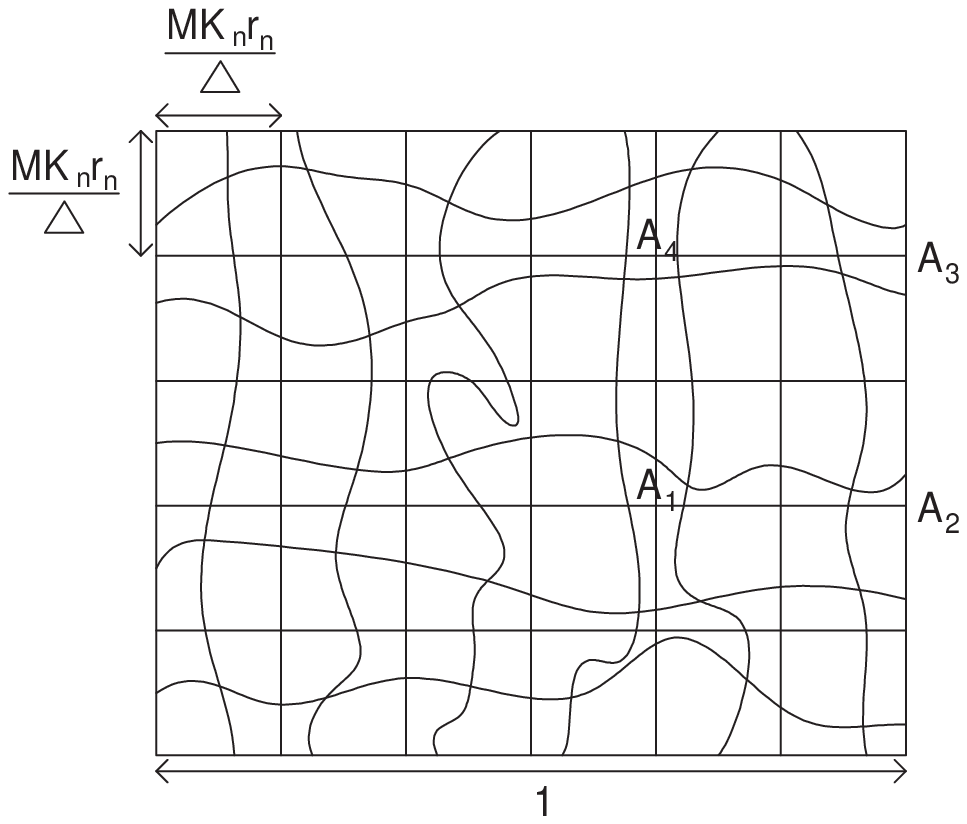}
   \nonumber
 }
\subfigure{
   \includegraphics[width = 2.5in] {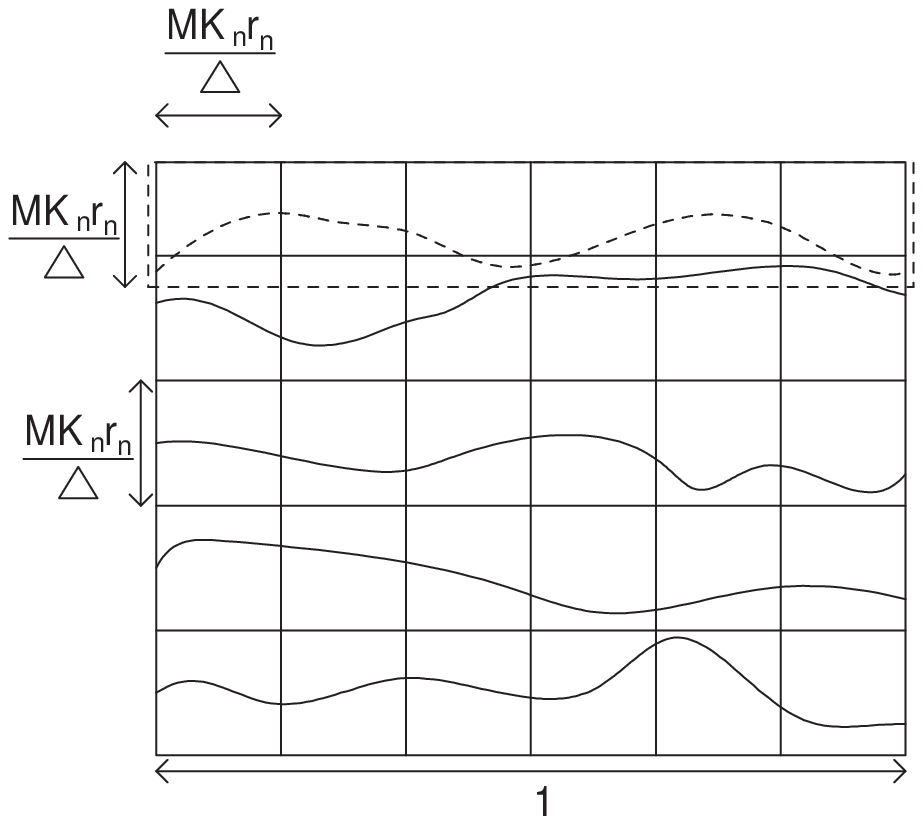}
   \nonumber
 }
\caption{Figures 3(a) and 3(b) of Ganesan (2013); replace \(\Delta\) with \(\epsilon_2^{-1}\) here. If \((MK_n\epsilon_2 r_n)^{-1}\) is not an integer, we start the tiling from the bottom and the two topmost rectangles overlap as above.}
\label{backbn}
\end{figure}

For convenience, we reproduce both the figures here in Figure~\ref{backbn}. We do an similar tiling for the vertical rectangles in \({\cal R}_V.\) If~\(R \in {\cal R}_H \cup {\cal R}_V,\) then~\(R\) contains exactly~\(MK_n\frac{1}{\epsilon_1 r_n}\) squares from~\(\{S_j\}\) and the total number of rectangles in~\({\cal R}_H \cup {\cal R}_V\) is
\begin{equation}\label{num_rect_est}
\#({\cal R}_H \cup {\cal R}_V) \leq 2[(MK_n\epsilon_1 r_n)^{-1}] + 2 \leq \frac{2}{MK_n\epsilon_1 r_n} + 2 \leq C_1\sqrt{n}
\end{equation}
for some constant \(C_1 > 0.\) As before \([x] \leq x\) is the largest integer less than or equal to \(x.\) The final estimate is obtained using the first inequality in (\ref{k_n_est}).

A square~\(S_j\) is said to be \emph{dense} if it contains at least~\(8\) vertices and sparse otherwise. A dense unoriented plus connected left right crossing is a set of distinct dense~\(\epsilon_1 r_n \times \epsilon_1 r_n\) squares~\((Y_1,\ldots,Y_D) \subseteq \{S_j\}\) contained in \(R\) satisfying the following properties.\\
(\(x1\)) The square~\(Y_1\) intersects the left side of~\(R\) and is plus adjacent (i.e., shares an edge) with~\(Y_2,\)\\
(\(x2\)) The square~\(Y_D\) intersects the right side of~\(R\) and is plus adjacent with~\(Y_{D-1}\) and\\
(\(x3\)) For every \(i, 2 \leq i \leq D-1,\) the square~\(Y_i\) is plus adjacent with~\(Y_{i-1}\) and~\(Y_{i+1}.\)\\
We have an analogous definition for star connected left right crossing by replacing plus adjacent above with star adjacent~(i.e., sharing a corner). We refer to Ganesan~(2015) for more on star and plus connected components.

For~\(R \in {\cal R}_H,\) let~\(F_n(R)\) be the event that the horizontally long rectangle~\(R \in {\cal R}_H\) contains an unoriented dense plus connected left right crossing of \(\epsilon_1r_n \times \epsilon_1r_n\) squares belonging to \(\{S_j\}.\) Analogously, for \(R \in {\cal R}_V,\) we define \(F_n(R)\) to be the event that \(R\) contains an unoriented plus connected occupied top bottom crossing. We have the following estimate. If \(R \in {\cal R}_H \cup {\cal R}_V,\) then
\begin{equation}\label{fn_r_est}
\mathbb{P}(F_n(R)) \geq 1 - \frac{1}{n^{10}}
\end{equation}
if \(M \geq 1\) is sufficiently large.\\
\emph{Proof of~(\ref{fn_r_est})}: If a plus connected dense top bottom crossing does not occur in~\(R,\) then there must exist a star connected sparse left right crossing. Fix any unoriented star connected left right crossing \(L_1 = (t_1,\ldots,t_l)\) containing \(l\) squares. Since the bottom edge of \(R\) has length \(MK_n\epsilon_1 r_n,\) we have that \(l \geq MK_n\) and since the left edge of \(R\) has length \(1\) there are \((\epsilon_1r_n)^{-1}\) possibilities for the square~\(t_1\) that starts from the left edge of \(R.\) For a fixed square \(t_1,\) there are at most~\(8^{l}\) choices for~\(L_1.\) For any fixed \(L_1,\) we have the following estimate
\begin{equation}\label{ti_vacant}
\mathbb{P}(t_i \text{ is sparse for }1 \leq i \leq l) \leq e^{-\theta l nr_n^2}
\end{equation}
for all \(n\) large and for some constant~\(\theta > 0.\)
To see~(\ref{ti_vacant}), we argue as follows. We have
\begin{equation}
\mathbb{P}\left(\bigcap_{i=1}^{l}\{t_i \text{ is sparse }\}\right) = \mathbb{P}\left(\bigcap_{j=1}^{n}\left\{X_j\notin \cup_{i} t_i \right\}\right) = \mathbb{P}\left(X_1 \notin \cup_i t_i\right)^n \label{a_sum5}
\end{equation}
where we recall that~\(X_i\) is the \(i^{th}\) random node placed in the unit square \(S\) according to the density \(f.\) In the above, \(\cup_i t_i\) is the union of the squares \(t_i\) and is a subset of the unit square \(S.\) The total area under the squares \(\cup_i t_i\) is~\(l \epsilon_1^2r_n^2\) and so we have
\begin{eqnarray}
\mathbb{P}\left(X_1 \notin \cup_i t_i\right)  = 1 - \int_{\cup_i t_i} f(x) dx \leq 1 -  l \epsilon_1^2r_n^2 \inf_{x \in S}f(x)  \leq 1- \theta l r_n^2,\nonumber
\end{eqnarray}
where \(\theta  := \frac{1}{16}\inf_{x \in S}f(x) > 0\) (see (\ref{f_eq})). The final inequality follows using \(\epsilon_1 \geq \frac{1}{4}\) (see the paragraph following (\ref{k_n_def})). Substituting into (\ref{a_sum5}) gives
\begin{equation}\label{vacant_est}
\mathbb{P}\left(t_i \text{ is sparse for } 1 \leq i \leq l\right) \leq (1- \theta  l r_n^2)^{n} \leq e^{-l\theta nr_n^2} \leq e^{-\theta l nr_n^2}
\end{equation}
where we use \(1-x \leq e^{-x}\) for all \(x > 0\) in obtaining the second estimate.

Let~\(LR_V\) denote the event that~\(R\) contains a star connected sparse left right crossing. Using~(\ref{ti_vacant}), we have
\begin{eqnarray}
\mathbb{P}(LR_V) &\leq& MK_n \sum_{l \geq MK_n} 8^{l} e^{-\theta l n r_n^2} \nonumber\\
&\leq& MK_n \frac{8e^{- \theta MK_n nr_n^2}}{1-8e^{-\theta nr_n^2}} \nonumber\\
&\leq& 16MK_ne^{-\theta M K_n nr_n^2} \nonumber
\end{eqnarray}
for all~\(n\) large. The final estimate is true since~\(nr_n^2 \rightarrow \infty\) and so~\(1-8e^{-\theta nr_n^2} \geq \frac{1}{2}\) for all~\(n\) large. Using~\(K_n nr_n^2 \geq \log{n}\) (see~(\ref{k_n_def})), we have that \[\mathbb{P}(LR_V) \leq 16MK_n e^{-\theta M\log{n}} \leq 32M\log{n}e^{-\theta M\log{n}} \leq \frac{1}{n^{9}}\] provided~\(M >0\) is large. The middle estimate above is true since~\(nr_n^2 \geq 1\) and so~\(K_n \leq 2\log{n}\) for all~\(n\) large. Since one of the events~\(F_n(R)\) or~\(LR_V\) must always occur, this proves~(\ref{fn_r_est}).~\(\qed\)

Fix~\(M \geq 1\) as in~(\ref{fn_r_est}) and set
\begin{equation}\label{q_n_def}
F_{n} := \bigcap_{R \in {\cal R}_H \cup {\cal R}_V} F_n(R).
\end{equation}
We have that
\begin{equation} \label{q_n_est}
\mathbb{P}(F_{n}) \geq 1 - \#({\cal R}_H \cup {\cal R}_V)\frac{1}{n^{10}} \geq 1 - C_1\sqrt{n}\frac{1}{n^{10}} \geq 1 - \frac{1}{n^9}
\end{equation}
for all~\(n\) large. The second inequality follows from~(\ref{num_rect_est}). We note that if~\(F_{n}\) occurs, we obtain a backbone of crossings containing vertices close to all sides of~\(S.\) In Figure~\ref{backbn}, the wavy lines correspond to the backbone. By considering lowermost occupied left right crossings of rectangles  in~\({\cal R}_H\) and leftmost top bottom crossings of rectangles in~\({\cal R}_V,\) we obtain a unique backbone of crossings which we call~\({\cal B}.\)

By construction, any two vertices in star adjacent dense squares of the backbone~\({\cal B}\) are connected by an edge and so the set of all vertices belonging to the squares in~\({\cal B}\) form a connected component of the graph~\(G\) which we denote by~\(C_G.\)


\subsection*{Estimating sizes of small components}
Let~\(C_G\) denote the component of the graph~\(G\) belonging to the backbone~\({\cal B}\) as defined in the previous subsection. Letting
\begin{equation}\label{xo_def}
X_O = \sum_{C \neq C_G} \#C = n-\#C_G
\end{equation}
denote the sum of sizes of all other components and arguing as in the proof of Lemma~\(3\) of Ganesan~(2013), we have that
\begin{equation}\label{xo_est}
\mathbb{E}X_O \leq ne^{-2\beta nr_n^2}
\end{equation}
for some constant~\(\beta > 0\) and for all~\(n\) large. We give a proof below for completeness.\\\\
\emph{Proof of~(\ref{xo_est})}: Let~\(M > 0\) be the constant as in the definition of the event~\(F_n\) (see~(\ref{q_n_def})). For~\(A \in \{S_j\},\) let~\(U_{2MK_n}(A)\) be the~\(2MK_n \epsilon_1 r_n \times 2MK_n \epsilon_1 r_n\) square with centre closest to the centre of~\(A\) and containing exactly~\((2MK_n)^2\) squares in~\(\{S_j\}.\) If there is more than one choice, we fix one according to a deterministic rule. For example, the centre with the least~\(x-\)coordinate and the least~\(y-\)coordinate.

By construction of the backbone, if~\(C \neq C_G\) is a component of the graph~\(G,\) then there is a square~\(A = A(C) \in \{S_j\}\) such that~\(C\) is contained in the bigger square~\(U_{2MK_n}(A).\) We therefore estimate the sizes of all components other than~\(C_G\) as follows. We define the random variable~\(X(A)\) as follows. Let~\(C(A)\) be the star connected dense component containing~\(A.\) We have that~\(C(A) = \emptyset\) if~\(A\) itself is sparse. If~\(A\) is dense, then let~\(V(A)\) be the event that every square in~\(C(A)\) is contained in the bigger square~\(W_{2MK_n}(A).\)

Letting~\(N(S_j)\) denote the number of vertices in the square~\(S_j \subset \{S_k\},\) define
\begin{equation}\label{xa_def}
X(A) = \sum_{S_j \in C(A)} N(S_j) \ind(V(A))
\end{equation}
to be the total number of vertices contained in the component~\(C(A).\) From the first statement in the previous paragraph, the term
\begin{equation}\label{yo_def}
Y_O = \sum_{A \in \{S_k\}} X(A)
\end{equation}
is an upper bound for~\(X_O\) defined in~(\ref{xo_def}). For a fixed~\(A \in \{S_j\},\) we have the following estimate.
\begin{equation}\label{xos_est}
\mathbb{E}X(A) \leq nr_n^2e^{-\beta_1 nr_n^2}
\end{equation}
for some constant~\(\beta_1 > 0\) and for all~\(n\) large.\\\\
\emph{Proof of~(\ref{xos_est})}: Suppose that the event~\(V(A)\) occurs and the component~\(C(A)\) contains~\(k\) squares. From the estimate~\((7)\) of Ganesan~(2013), we have
\begin{equation}\label{ca_k_est}
\mathbb{P}\left(\{\#C(A) = k\} \bigcap V(A)\right) \leq ke^{-\theta_1 nr_n^2 \sqrt{k}}
\end{equation}
for some constant~\(\theta_1 > 0\) and for all~\(n \geq N_1.\) Here~\(\theta_1\) and~\(N_1\) do not depend on~\(k.\) Proceeding as in the analysis following~\((7)\) of Ganesan~(2013), we then obtain~(\ref{xos_est}).

For completeness we give a proof of~(\ref{ca_k_est}).\\
\emph{Proof of~(\ref{ca_k_est})}: We write
Suppose \(C(A)\) contains \(k\) squares. We use Theorem~\(1\) of Ganesan (2015) and obtain that the outermost boundary~\(\partial_A\) of \(C(A)\) is a connected union of cycles \(\cup_{i = 1}^{h} H_i\) each consisting only of boundary edges; i.e., edges either contained in the boundary of the unit square~\(S\) or edges present in the interior of~\(S\) and adjacent to one sparse and one dense square of~\(C(A).\) Moreover, there is a circuit~\(\Pi\) consisting of the edges of~\(C(A).\) By a circuit of edges, we mean a sequence of distinct edges \((e_1,\ldots,e_k)\) such that the following three statements hold: The edge \(e_i\) shares one endvertex with \(e_{i+1}\) and one endvertex with \(e_{i-1}\) for all \(2 \leq i \leq k-1.\) The edge \(e_k\) shares one endvertex with \(e_1\) and one endvertex with \(e_{k-1}\) and the edge \(e_1\) shares one endvertex with \(e_k\) and one endvertex with \(e_2.\)


We have the following properties regarding the circuit~\(\Pi.\)\\
\((l1)\) If \(\#\Pi\) denote the (random) number of edges in \(\Pi,\) we have that  \(\frac{\sqrt{k}}{4} \leq \#\Pi \leq 4k.\)\\
\((l2)\) If \(N_{vac}\) denotes the number of distinct sparse squares sharing a edge with some occupied square of \(C(A),\) we have that \(N_{vac} \geq \frac{\#\Pi}{8}.\)\\
\((l3)\) Every edge in~\(\Pi\) is contained in the larger square~\(U_{2MK_n}(A).\)\\

\emph{Proof of \((l1)-(l3)\)}: The property~\((l3)\) is true by definition. For the upper bound in property \((l1)\), we use the fact that each occupied square contains four edges and every edge in \(\Pi\) is adjacent to some occupied square of~\(C(A).\)  To see the lower bound, we suppose \(n_i\) squares of \(C(A)\) is contained in the interior of the cycle \(H_i, 1 \leq i \leq h\) so that \(\sum_{i=1}^{h} n_i = k.\) The final equality is true since every square of \(C(A)\) is contained in the interior of some cycle, by Theorem~1 of Ganesan (2015).

If cycle \(H_i\) contains \(l_i\) edges, we must then have that \(n_i \leq 16l_i^2.\) To see this fix any vertex \(v \in H_i.\) All the \(\epsilon_1 r_n \times \epsilon_1 r_n\) squares contained in the interior of~\(H_i\) are contained in interior of the bigger \(4l_i\epsilon_1 r_n \times 4l_i \epsilon_1 r_n\) square centred at~\(v.\) Therefore the total number of squares \(n_i\) is at most \(16l_i^2.\) Summing over~\(i\) gives \[k = \sum_{i=1}^{h} n_i \leq 16 \sum_{i=1}^{h}l_i^2 \leq 16 \left(\sum_{i=1}^{h} l_i\right)^2 = 16\left(\#\Pi\right)^2.\] This implies that \(\#\Pi \geq \frac{\sqrt{k}}{4}.\)

To see \((l2)\) is true, we suppose as above that the cycle \(H_i\) of the outermost boundary contains \(l_i\) edges, \(1 \leq i \leq h.\) We write \(l_i = l_{i,1} + l_{i,2},\) where \(l_{i,1}\) is the number of edges of \(H_i\) contained in the boundary of the unit square~\(S.\) Suppose edge \(e \in H_i\) touches the left edge of~\(S\) and suppose \((x_e,y_e)\) is the centre of the segment formed by the edge \(e.\) The cycle \(H_i\) cuts the line \(y = y_e\) at some unique edge \(e_1 = e_1(e)\) contained in the interior of \(R^{left}_T.\) Thus \(l_{i,2} \geq l_{i,1}\) and so \(2 l_{i,2} \geq l_i.\) Summing over \(i\) gives \(2\sum_{i=1}^{h} l_{i,2} \geq \#\Pi.\)

The term \(\sum_{i=1}^{h} l_{i,2}\) denotes the number of edges of \(\Pi\) contained in the interior of the unit square \(S.\) Every such edge obtained is adjacent to a sparse square and an occupied square of~\(C(A).\) Since each sparse square has four edges, the number of \emph{distinct} sparse squares attached to some edge of~\(\Pi\) is~\(N_{vac} \geq \frac{\#\Pi}{8}.\)~\(\qed\)

Using properties \((l1)-(l2)\) we have for a fixed \(k \geq 1\) that
\begin{equation}
\mathbb{P}\left(\{\#C(A) =k\} \cap V(A)\right) = \sum_{\pi : \frac{\sqrt{k}}{4} \leq \#\pi \leq 4k } \mathbb{P}\left(\{\#C(A) =k\} \cap V(A) \cap \{\Pi = \pi\}\right)\label{a_sum1}
\end{equation}
where the summation is over all circuits \(\pi\) surrounding the square~\(A,\) and contained in the larger square~\(U_{2MK_n}(A)\) (see property~\((l3).\) For a realization~\(\Pi = \pi\) with \(\#\pi =l,\) the set of sparse \(\epsilon_1 r_n \times \epsilon_1 r_n\) squares containing  some edge of~\(\pi\) and lying in the exterior of every cycle of~\(\pi,\) is fixed. Letting \(t_i = t_i(\pi),1 \leq i \leq n_{vac}\) be the set of such sparse squares, we have
\begin{equation}
\mathbb{P}\left(\{\#C(A) =k\} \cap V(A) \cap \{\Pi = \pi\}\right) \leq \mathbb{P}\left(\bigcap_{i=1}^{n_{vac}}\{t_i \text{ is sparse }\} \right). \label{a_sum3}
\end{equation}

Using the estimate~(\ref{ti_vacant}) in (\ref{a_sum1}) gives
\begin{eqnarray}
\mathbb{P}\left(\{\#C(A) =k\} \cap V(A) \right) &\leq& \sum_{\frac{\sqrt{k}}{4} \leq l \leq 4k} \sum_{\pi : \#\pi = l}e^{-\theta_3 l nr_n^2} \nonumber\\
&\leq& \sum_{\frac{\sqrt{k}}{4} \leq l \leq 4k} l.8^{l} e^{-\theta l nr_n^2} \nonumber\\
&\leq& 4k \sum_{\frac{\sqrt{k}}{4} \leq l \leq 4k} 8^{l} e^{-\theta l nr_n^2} \label{count_pi}
\end{eqnarray}
The middle inequality is obtained using the fact that the number of circuits of length~\(l\) surrouding~\(A\) is at most \(l.8^{l}.\) To see this is true, we draw axes parallel to the sides of~\(A\) such that one corner of the square~\(A\) is the origin. The circuit~\(\pi\) intersects the~\(X-\)axis at some point~\(g(\pi).\) The number of choices for~\(g(\pi)\) is at most~\(l\) and for each fixed choice of~\(g(\pi),\) the number of choices for~\(\pi\) is at most~\(8^{l}.\) This proves~(\ref{count_pi}).

From~(\ref{count_pi}), we have \[\mathbb{P}\left(\{\#C(A) =k\} \cap V(A )\right) \leq \frac{4k}{1-8e^{-\theta nr_n^2}}\left(8.e^{-\theta nr_n^2}\right)^{\sqrt{k}/4} \leq k.e^{-\theta_1 nr_n^2 \sqrt{k}}\] for all~\(n \geq N_1.\) Here~\(0 < \theta_1 < \theta\) is fixed and~\(N_1 \geq 1\) does not depend on~\(k.\) The final estimate is obtained using~\(nr_n^2 \rightarrow \infty\) as~\(n \rightarrow \infty.\)~\(\qed\)

By our choice of~\(\epsilon_1 \in \left(\frac{1}{4},\frac{1}{5}\right)\) in the first paragraph of this proof, the number of squares in~\(\{S_j\}\) is~\(\frac{1}{\epsilon_1^2 r_n^2} \leq \frac{25}{r_n^2}\) and so we have from~(\ref{yo_def}) that \[\mathbb{E}Y_O \leq 25 n e^{-\beta_1 nr_n^2} \leq ne^{-\beta_2 nr_n^2}\] for a fixed constant~\(0 < \beta_2 < \beta_1\) and for all~\(n\) large. The final inequality is true since~\(nr_n^2 \rightarrow \infty\) as~\(n \rightarrow \infty\). Since~\(X_O \leq Y_O,\) we obtain~(\ref{xo_est}).~\(\qed\)

Using Markov inequality and~(\ref{xo_def}), we have
\begin{equation}\label{cg_est}
\mathbb{P}\left(X_O \geq ne^{-\beta nr_n^2}\right) \leq e^{-\beta nr_n^2}
\end{equation}
for all~\(n\) large. Let~\(F_n\) be the event (see~(\ref{q_n_def})) that a backbone of dense crossings occur in the unit square~\(S.\) Defining the event
\begin{equation}\label{en_def}
E_n = F_n \bigcap \{X_O \leq ne^{-\beta nr_n^2}\}
\end{equation}
we have from~(\ref{q_n_est}) that
\begin{equation}\label{en_est}
\mathbb{P}(E_n) \geq 1-\frac{1}{n^9}-e^{-\beta nr_n^2} \geq 1-e^{-\beta_1 nr_n^2}
\end{equation}
for all~\(n\) large and for some constant~\(\beta_1 > 0.\) If~\(E_n\) occurs, then the component~\(C_G\) belonging to the backbone~\({\cal B}\) contains at least~\(n-ne^{-\beta nr_n^2}\) vertices and is therefore the largest component.

\subsection*{Forming the long cycle using the backbone}
Suppose that the event~\(E_n\) defined in~(\ref{en_def}) occurs. From~(\ref{en_est}), we have that the largest component~\(C_G\) contains at least~\(n-ne^{-\beta nr_n^2}\) vertices with probability at least~\(1-e^{-\beta_1 nr_n^2}.\) Moreover every node of the component~\(C_G\) belongs to some dense square in the backbone~\({\cal B}.\) Letting~\({\cal B} = \{Y_i, 1 \leq i \leq t\},\) we inductively obtain a cycle of edges  in the graph~\(G\) containing all vertices of~\(C_G.\)

Consider a sequence of star connected components~\({\cal B}_i, 1 \leq i \leq t\) such that~\({\cal B}_1 = \{W_1\} = \{Y_1\}, {\cal B}_t = {\cal B}\) and for~\(1 \leq i \leq t-1,\) the component~\({\cal B}_{i+1}\) contains one more square~\(W_{i+1} \subset \{Y_j\}\) than~\({\cal B}_i.\) The square~\(W_{i+1}\) is star adjacent to some square~\(W_l \in {\cal B}_i.\) Thus~\({\cal B}_i= \cup_{1 \leq j \leq i} W_j\) and for~\(1 \leq i \leq t,\) let~\(\eta_i\) be a cycle containing all the vertices present in the square~\(W_i.\)

We set~\(\tau_1 = \eta_1\) and iteratively construct cycles~\(\tau_i,1 \leq i \leq t,\) using~\(\{\eta_i\}_{1 \leq i \leq t}.\) The final cycle~\(\tau_t\) is then the desired long cycle. For~\(1 \leq i \leq t,\) we have the following properties for the cycle~\(\tau_i = (g_1,\ldots,g_w)\) where each~\(g_i\) is an edge.\\
\((b1)\) All edges of~\(\{\eta_j\}_{1 \leq j \leq i}\) not removed so far in the iteration process belong to the cycle~\(\tau_{i}.\)\\
\((b1)\) Let~\(W_{i+1} = {\cal B}_{i+1} \setminus {\cal B}_i\) be adjacent to some square~\(W_{l} \in {\cal B}_i.\) Here~\(1 \leq l \leq i\) and there exists an edge~\(g_l \in \eta_l \cap \tau_i.\)\\
\emph{Proof of~\((b1)-(b2)\) for~\(i = 1\)}: The square~\(W_l\) contains at least~\(8\) vertices and so~\((b1)-(b2)\) is true.\(\qed\)

Using properties~\((b1)-(b2),\) we form the new cycle~\(\tau_{i+1}\) as follows. Let~\(u_l\) and~\(v_l\) be the endvertices of edge~\(g_l\) which belong to the dense square~\(W_l.\) We recall that~\(\eta_{i+1}\) is a cycle of edges containing all the vertices in the square~\(W_{i+1}.\) Remove one edge from~\(\eta_{i+1}\) and let~\(a\) and~\(b\) be the endvertices of resulting path~\(P_{i+1}.\)

\begin{figure}[tbp]
\centering
\includegraphics[width=2.5in, trim= 130 290 200 230, clip=true]{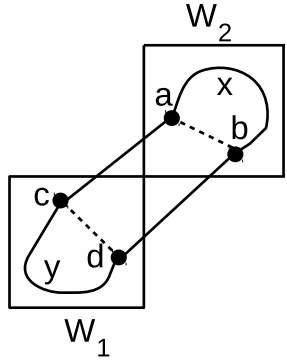}
\caption{Merging the cycle~\(\tau_1 = cydc\) contained in the square~\(W_1\) and the cycle~\(\eta_2 = axba\) contained in the square~\(W_2.\)}
\label{fig_long_rgg}
\end{figure}

The vertices~\(a\) and~\(u_l\) belong to star adjacent squares in~\(\{S_j\}\) and are therefore connected by an edge. Similarly the vertices~\(b\) and~\(v_l\) are also connected by an edge. We then merge the path~\(\tau_i\setminus\{g_l\}\) with the path~\(P_{i+1}\) to get the new cycle
\begin{equation}\label{di_1_cyc}
\tau_{i+1} = \left(\tau_i \setminus \{g_l\}\right) \cup P_{i+1}.
\end{equation}
This is illustrated for~\(i =l = 1\) in Figure~\ref{fig_long_rgg}, where the cycle~\(\tau_1 = \eta_1\) contained in the square~\(W_1\) is given by the wavy path~\(cdyc\) with~\(cd\) denoting the edge~\(g_1.\) The cycle~\(\eta_2 = axba\) and the path~\(P_2 = axb.\) The new cycle~\(\tau_2 = cydbxac.\)


The cycle~\(\tau_{2}\) contains all the vertices in the component~\({\cal B}_2\) and also satisfies properties~\((b1)-(b2).\) We continue this process iteratively and the cycle~\(\tau_i\) obtained at the end of iteration~\(i\) also satisfies properties~\((b1)-(b2).\)\\
\emph{Proof of~\((b1)-(b2)\) for~\(i \geq 2\)}: The proof of~\((b1)\) is true by construction. To prove~\((b2),\) we argue as follows. At the end of each iteration at most one edge is removed from each cycle~\(\eta_j, 1 \leq j \leq t.\) The square~\(W_l\) star adjacent to~\(W_{i+1}\) contains at least~\(8\) vertices and so the corresponding cycle~\(\eta_l\) containing all the vertices of~\(W_l\) has at least~\(8\) edges. There are exactly~\(8\) squares star adjacent to~\(W_l\) and since~\(W_{i+1} \in {\cal B}_{i+1} \setminus {\cal B}_i\) is also star adjancent to~\(W_l,\) at most~\(7\) squares in~\({\cal B}_i\) are star adjacent to~\(W_l.\) This means that at most~\(7\) edges from~\(\eta_l\) have been removed so far in the iterative process above.~\(\qed\)

\section{Proof of~(\ref{eq_lg2_rgg}) in Theorem~\ref{long_rgg}}\label{pf_long_cycle_rgg_ham}
The proof is analogous as in the previous case with some minor modifications.
Suppose~\[nr_n^2 = \log{n} + 7\log{\log{n}} + \omega_n\] where~\(\omega_n \rightarrow \infty\) and~\(\frac{\omega_n}{\log{\log{n}}} \rightarrow 0\) as~\(n \rightarrow \infty.\) 

Divide~\(S\) into squares~\(\{S_j\}\) of side length~\(t_n,\) where
\begin{equation}\label{tn_def}
8nt^2_n =\log{n} + 7\log{\log{n}} + \omega_n - \delta_n
\end{equation}
and~\(\delta_n \in(1,2)\) is such that~\(\frac{1}{t_n}\) is an integer. The number~\(t_n\) is slightly less than~\(\frac{r_n}{2\sqrt{2}}\) and so if squares~\(S_{j_1}\) and~\(S_{j_2}\) share a corner, then every vertex in~\(S_{j_1}\) is joined to every vertex in~\(S_{j_2}\) by an edge. For a fixed square~\(S_j,\) we say that~\(S_j\) is dense if it contains at least~\(8\) vertices and sparse otherwise. Let~\(E(j)\) be the event that~\(S_j\) is sparse. If~\(S_{j_1},\ldots,S_{j_q}\) are fixed squares,~\(q\) not depending on~\(n,\) then
\begin{equation}\label{sj_sparse}
\mathbb{P}\left(\bigcap_{i=1}^{q}E(j_i)\right) \leq \frac{C(\log{n})^{8-\frac{7q}{8}}}{n^{q/8}}\exp\left(-\frac{q\omega_n}{8}\right)
\end{equation}
for some constant~\(C  = C(q) > 0\) and for all~\(n\) large.\\
\emph{Proof of~(\ref{sj_sparse})}: We have
\begin{eqnarray}
\mathbb{P}\left(\bigcap_{i=1}^{q}E(j_i)\right) &=& \sum_{k=0}^{8} {n \choose k} (qt_n^2)^{k} (1-qt_n^2)^{n-k} \nonumber\\
&\leq& \sum_{k=0}^{8} (qnt_n^2)^{k} (1-qt_n^2)^{n-k} \nonumber\\
&\leq& \frac{1}{(1-qt_n^2)^{8}}\sum_{k=0}^{8} (qnt_n^2)^{k} (1-qt_n^2)^{n} \nonumber\\
&\leq& \frac{1}{(1-qt_n^2)^{8}}\sum_{k=0}^{8} (qnt_n^2)^{k} e^{-qnt_n^2}.\label{sparse_est1}
\end{eqnarray}
The first inequality is obtained using~\({n \choose k} \leq n^{k}\) and the final inequality is obtained using the inequality~\(1-x < e^{-x}\) and the fact that~\(t_n < 1\) for all~\(n\) large~(see~(\ref{tn_def})).

From~(\ref{tn_def}), we in fact have that~\(t_n \rightarrow 0\) as~\(n \rightarrow \infty\) and so for a fixed~\(q\) not depending on~\(n,\) we have that~\((1-qt_n^2)^{-8} \leq 2\) for all~\(n\) large and so we have from~(\ref{sparse_est1}) that
\begin{equation}\label{sparse_est2}
\mathbb{P}\left(\bigcap_{i=1}^{q}E(j_i)\right) \leq 2\sum_{k=0}^{8} (qnt_n^2)^{k} e^{-qnt_n^2}.
\end{equation}
From the definition~(\ref{tn_def}), we also have that~\(qnt_n^2 \leq q\log{n}\) for all~\(n\) large and so from~(\ref{sparse_est2}) we have
\begin{eqnarray}
\mathbb{P}(E_j) &\leq& 18(q\log{n})^{8} e^{-qnt_n^2} \nonumber\\
&=& \frac{18q^8(\log{n})^{8}}{n^{q/8}(\log{n})^{7q/8}}\exp\left(-\frac{q\omega_n}{8} + \frac{q\delta_n}{8}\right) \nonumber\\
&\leq& \frac{C_1(\log{n})^{8-\frac{7q}{8}}}{n^{q/8}}\exp\left(-\frac{q\omega_n}{8}\right)\label{sparse_est3}
\end{eqnarray}
for all~\(n\) large and for some constant~\(C_1 > 0.\) The middle equality is obtained by substituting the expression for~\(nt_n^2\) from~(\ref{tn_def}) and the final estimate is obtained using the fact that~\(\delta_n \in (1,2).\)~\(\qed\)

\subsection*{Constructing the backbone}
As in the previous section, divide~\(S\) into disjoint~\(1 \times Mt_n\) horizontal rectangles and call the resulting set of rectangles as~\({\cal R}_H.\) Similarly divide~\(S\) into~\(Mt_n \times 1\) vertical rectangles and call the resulting set~\({\cal R}_V.\) Assume that the tiling of~\(S\) into rectangles in~\({\cal R}_H \cup {\cal R}_V\) is either as in Figure~\ref{backbn}(a) or as in Figure~\ref{backbn}(b) so that the number of rectangles in~\({\cal R}_H \cup {\cal R}_V\) is at most
\begin{equation}\label{num_rect_est2}
\frac{2}{Mt_n} + 2 \leq \frac{2}{M} \sqrt{\frac{8n}{\log{n}}} + 2 \leq \sqrt{n}
\end{equation}
for all~\(n\) large. The middle inequality is obtained using~\(8nt_n^2 \geq \log{n}\)~(see~(\ref{tn_def})).

For~\(R \in {\cal R}_H,\) let~\(F_n(R)\) be the event that the horizontally long rectangle~\(R \in {\cal R}_H\) contains an unoriented dense plus connected left right crossing of \(t_n \times t_n\) squares belonging to \(\{S_j\}.\) Analogously, for \(R \in {\cal R}_V,\) we define \(F_n(R)\) to be the event that \(R\) contains an unoriented plus connected occupied top bottom crossing. Analogous to the proof of~(\ref{fn_r_est}), we have that if \(R \in {\cal R}_H \cup {\cal R}_V,\) then \[\mathbb{P}(F_n(R)) \geq 1 - \frac{1}{n^{10}}\] if \(M \geq 1\) is a constant sufficiently large. Fixing such an \(M\) and setting
\begin{equation}\label{q_n_def2}
F_{n} := \bigcap_{R \in {\cal R}_H \cup {\cal R}_V} F_n(R),
\end{equation}
we have that
\begin{equation} \label{q_n_est2}
\mathbb{P}(F_{n}) \geq 1 - \#({\cal R}_H \cup {\cal R}_V)\frac{1}{n^{10}} \geq 1 - \sqrt{n}\frac{1}{n^{10}} \geq 1 - \frac{1}{n^9}
\end{equation}
for all \(n\) large. The second inequality follows from (\ref{num_rect_est2}). We note that if \(F_{n}\) occurs, we obtain a backbone of crossings containing vertices close to all sides of \(S.\) In Figure~\ref{backbn}, the wavy lines correspond to the backbone. By considering lowermost occupied left right crossings of rectangles  in \({\cal R}_H\) and leftmost top bottom crossings of rectangles in \({\cal R}_V,\) we obtain a unique backbone of crossings which we call~\({\cal B}.\)

\subsection*{Finding isolated dense components}
Suppose that the event~\(F_n\) defined in~(\ref{q_n_def2}) occurs and let~\({\cal B}\) be the corresponding backbone constructed above. For a square~\(A \in \{S_j\},\) let~\(C(A)\) be the star connected dense component containing~\(A.\) Define
\begin{equation}\label{za_def}
I(A) = F_n \bigcap \{C(A) \neq {\cal B}\} \bigcap \{A \text{ is dense}\}
\end{equation}
to be the event that the dense component containing~\(A\) is not the backbone~\({\cal B}.\) The existence of the backbone is guaranteed by the occurrence of the event~\(F_n.\) Let
\begin{equation}\label{ei_def}
I_n = \bigcup_{A \in \{S_j\}} I(A)
\end{equation}
be the event that there exists a dense component that is not equal to the backbone~\({\cal B}.\) 


We have
\begin{equation}\label{ei_main_est}
\mathbb{P}(I_n) \leq Ce^{-\omega_n}
\end{equation}
for some constant~\(C > 0\) and for all~\(n\) large. Here~\(\omega_n \longrightarrow \infty\) is as in~(\ref{rn_hom}).\\\\
\emph{Proof of~(\ref{ei_main_est})}: We evaluate the probability of the event~\(I_n\) by estimating the size of each isolated dense component~\(C(A).\) We have some notations first. Let~\(L = (A_1,\ldots,A_t)\) be a sequence of distinct squares in~\(\{S_j\}.\) We say that~\(L\) is a \emph{plus connected~\(S-\)cycle} if the following conditions are satisfied.
\((a)\) For~\(1 \leq i \leq t-1,\) the square~\(A_i\) is plus adjacent  (i.e., shares an edge) with the square~\(A_{i+1}.\)\\
\((b)\) The square~\(A_t\) is plus adjacent to~\(A_{t-1}\) and~\(A_1.\)

Let~\(S(1+2t_n)\) be the larger square with same centre as the unit square~\(S\) and of side length~\(1+2t_n.\) The set of squares obtained by tiling~\(S(1+2t_n)\) into~\(t_n \times t_n\) squares is~\(\{S_j\} \cup\{Q_j\}_{j=1}^{w}\) where~\((Q_1,\ldots,Q_w)\) is an~\(S-\)cycle of squares lying in the exterior of~\(S\) and intersecting~\(S.\) We define every square in~\(\{Q_j\}\) to be sparse. 

Fix~\(A \in \{S_j\}\) and suppose that the event~\(I(A)\) occurs; i.e., there is a backbone~\({\cal B}\) containing dense squares and the dense component~\(C(A)\) containing~\(A\) is not~\({\cal B}.\) By construction of the backbone~\({\cal B},\) every square in~\(C(A)\) is contained in the~\(2Mt_n \times 2Mt_n\) bigger square~\(U_{2M}(A).\) Here~\(M \geq 1\) is the constant in~(\ref{q_n_def2}) and as in the proof of~(\ref{xo_est}), the square~\(U_{2M}(A)\) is the~\(2M t_n \times 2Mt_n\) square with centre closest to the centre of~\(A\) and containing exactly~\((2M)^2\) squares in~\(\{S_j\}.\)

From Theorem~\(1\) of Ganesan~(2015), we have that there is a plus connected~\(S-\)cycle~\(L_{cyc} = (R_1,\ldots,R_T)\) of sparse squares in~\(\{S_j\}\cup\{Q_j\}\) surrounding~\(C(A).\) Since~\(C(A)\) is contained in~\(U_{2M}(A)\) (see previous paragraph), we have the following property.
\begin{equation}\label{lcyc}
\text{Every square of~\(L_{cyc}\) is contained in the~\(3Mt_n \times 3Mt_n\) square~\(U_{3M}(A).\)}
\end{equation}
We consider three cases below depending on where the square~\(A\) is located. \\\\
\underline{\emph{Case I}}: The square~\(A\) intersects one of the corners of the unit square~\(S.\) Fix a realization~\(L_{cyc} = \pi\) where~\(\pi = (S_{j_1},\ldots,S_{j_w})\) is a deterministic~\(S-\)cycle surrounding the square~\(A\) and contained in the bigger~\(3Mt_n \times 3Mt_n\) square\\\(U_{3M}(A).\) Let~\({\cal T}_A\) denote the set of all such plus connected~\(S-\)cycles. We have that
\begin{equation}\label{n_pi_est}
\#{\cal T}_A \leq \sum_{w=1}^{(3M)^2}w.8^w
\end{equation}
\emph{Proof of~(\ref{n_pi_est})}: Fix a~\(\pi \in {\cal T}_A\) containing~\(w\) squares. Consider axes parallel to the sides of the square~\(A\) with origin denoted by one of the corners of~\(A.\) Let~\(g(\pi) \subset\{S_j\}\) be the square intersecting the positive~\(X-\)axis. There are at most~\(w\) choices for~\(g(\pi)\) since there are at most~\(w\) squares in~\(\pi.\) For each fixed choice of~\(g(\pi),\) there are at most~\(8^{w}\) choices for the~\(S-\)cycle~\(\pi.\) Thus there are at most \(w.8^{w}\) choices for the cycle~\(\pi.\) 

From~(\ref{lcyc}), we also have that every square in~\(\pi\) is contained within the~\(3Mt_n \times 3Mt_n\) square~\(U_{3M}(A)\) and so~\(w \leq (3M)^2.\) This proves~(\ref{n_pi_est}).~\(\qed\)

Every square in the \(S-\)cycle~\(\pi\) is sparse and at least three squares of~\(\pi\) must lie in the interior of the unit square~\(S.\) Let~\(S_{i_1},S_{i_2}\) and~\(S_{i_3}\) be the sparse squares with least such indices. Recalling from~(\ref{sj_sparse}) that~\(E(i_j)\) denotes the event that~\(S_{i_j}\) is sparse, we have
\begin{eqnarray}
\mathbb{P}\left(I(A)\right) &\leq& \sum_{\pi \in {\cal T}_A} \mathbb{P}\left(\left\{L_{cyc}  = \pi\right\} \bigcap \bigcap_{j=1}^{3}E(i_j)\right) \nonumber\\
&\leq& \sum_{\pi \in {\cal T}_A} \mathbb{P}\left(\bigcap_{j=1}^{3}E(i_j)\right) \nonumber\\
&\leq& \sum_{\pi \in {\cal T}_A} D_0\frac{(\log{n})^{43/8}}{n^{3/8}}\exp\left(-\frac{3\omega_n}{8}\right) \label{pza1_0}\\
&\leq& D\frac{(\log{n})^{43/8}}{n^{3/8}}\exp\left(-\frac{3\omega_n}{8}\right).\label{pza1_1}
\end{eqnarray}
The estimate~(\ref{pza1_0}) follows from~(\ref{sj_sparse}) by setting~\(q = 3\) and~\(D_0 > 0\) is the constant in~(\ref{sj_sparse}). In~(\ref{pza1_1}), the constant~\(D =  \sum_{w=1}^{(3M)^2}w.8^wD_0\) and the estimate~(\ref{pza1_1}) follows from~(\ref{n_pi_est}).\\\\
\underline{\emph{Case II}}: The square~\(A\) does not intersect any corner of the unit square~\(S\) but is within a distance of~\(3t_n\) from the boundary of~\(S.\)\\


In this case at least~\(5\) squares in the \(S-\)cycle~\(L_{cyc}\) lie in the interior of the unit square~\(S.\) Arguing as in Case~\((I)\) above and using~(\ref{sj_sparse}) with~\(q = 5,\) we have
\begin{equation}\label{pza2_2}
\mathbb{P}\left(I(A)\right) \leq D \frac{(\log{n})^{29/8}}{n^{5/8}}\exp\left(-\frac{5\omega_n}{8}\right)
\end{equation}
for some constant~\(D > 0.\)\\\\
\underline{\emph{Case III}}: The square~\(A\) is at a distance of~\(3t_n\) away from the boundary of~\(S.\)


In this case at least~\(8\) squares in the \(S-\)cycle~\(L_{cyc}\) lie in the interior of the unit square~\(S.\) Arguing as in Case~\((I)\) above and using~(\ref{sj_sparse}) with~\(q = 8,\) we have
\begin{equation}\label{pza3_2}
\mathbb{P}\left(I(A)\right) \leq  D\frac{\log{n}}{n}e^{-\omega_n}
\end{equation}
for some constant~\(D > 0.\)\\

Let~\(N_j\) be the number of squares satisfying Case~\((j)\) for~\(j\in \{I, II,III\}.\) We have that
\begin{equation}\label{nj_est}
N_I = 4, N_{II} \leq \sqrt{n} \text{ and }N_{III} \leq \frac{8n}{\log{n}}.
\end{equation}
\emph{Proof of~(\ref{nj_est})}: The first estimate on~\(N_I\) is true since there are four corners of~\(S.\) For~\(N_{II},\) we have that the number of squares intersecting the boundary of~\(S\) and contained in the interior of~\(S\) is at most~\(\frac{4}{t_n}.\) Therefore the number of squares at a distance of at most~\(3t_n\) from the boundary of~\(S\) is at most~\[\frac{12}{t_n} \leq 12\sqrt{\frac{8n}{\log{n}}} \leq \sqrt{n}\] for all~\(n\) large. The middle inequality follows since~\(8nt_n^2 \geq \log{n}\) for all~\(n\) large (see~(\ref{tn_def})).

Similarly, the final estimate on~\(N_{III}\) is true since the total number of squares in~\(\{S_j\}\) contained in the interior of~\(S\) is~\(\frac{1}{t_n^2} \leq \frac{8n}{\log{n}}.\)~\(\qed\) 

Using~(\ref{nj_est}), we have from~(\ref{ei_def}),~(\ref{pza1_1}),~(\ref{pza2_2}) and~(\ref{pza3_2}) that
\begin{eqnarray}
\mathbb{P}\left(I_n \right) &\leq& D \frac{(\log{n})^{43/8}}{n^{3/8}}\exp\left(-\frac{3\omega_n}{8}\right) + \sqrt{n}D \frac{(\log{n})^{29/8}}{n^{5/8}}\exp\left(-\frac{5\omega_n}{8}\right) \nonumber\\
&&\;\;\;\;\;\;\;\;\;\;+\;\;\;\frac{8n}{\log{n}} D\frac{\log{n}}{n}e^{-\omega_n}. \label{ei_esta}
\end{eqnarray}
We have that
\begin{equation}\label{ei_est1}
\frac{(\log{n})^{43/8}}{n^{3/8}}\exp\left(-\frac{3\omega_n}{8}\right) \leq e^{-\omega_n}
\end{equation}
and
\begin{equation}\label{ei_est2}
\sqrt{n}\frac{(\log{n})^{29/8}}{n^{5/8}}\exp\left(-\frac{5\omega_n}{8}\right)  \leq e^{-\omega_n}
\end{equation}
for all~\(n\) large.\\
\emph{Proof of~(\ref{ei_est1}) and~(\ref{ei_est2})}: We prove~(\ref{ei_est1}) and the proof of~(\ref{ei_est2}) is analogous.
To prove~(\ref{ei_est1}), it is enough to see that
\[\exp\left(\frac{5\omega_n}{8}\right) \leq \frac{n^{3/8}}{(\log{n})^{43/8}}\] or equivalently that
\[\frac{5\omega_n}{8} \leq \frac{3}{8} \log{n} - \frac{43}{8}\log{\log{n}}\] which is true for all~\(n\) large since~\(\frac{\omega_n}{\log{\log{n}}} \longrightarrow 0\) as~\(n \rightarrow \infty.\)~\(\qed\)

Using~(\ref{ei_est1}) and~(\ref{ei_est2}) into~(\ref{ei_esta}) gives~(\ref{ei_main_est}).~\(\qed\)


\subsection*{Isolated sparse squares}
Let~\(A \in \{S_j\}\) be any square and let~\(J(A)\) be the event that all the squares star adjacent to~\(A\) and contained in the unit square~\(S\) are sparse. Defining
\begin{equation}\label{fi_def}
J_n = \bigcup_{A \in \{S_j\}} J(A)
\end{equation}
we have that
\begin{equation}\label{fi_main_est}
\mathbb{P}(J_n) \leq C e^{-\omega_n}
\end{equation}
for some constant~\(C > 0\) and for all~\(n\) large. In particular if the event~\(J_n^c\) occurs, then every sparse square is star adjacent to some dense square.\\\\
\emph{Proof of~(\ref{fi_main_est})}: To estimate the probability of the event~\(J_n,\) we consider cases~\(I,II\) and~\(III\) as in the previous subsection. In case~\(I\) there are three squares star adjacent to~\(A\) and contained in the unit square. Using~(\ref{sj_sparse}) with~\(q = 3,\) we therefore have
\begin{equation}\label{case1}
\mathbb{P}(J(A)) \leq D_0 \frac{(\log{n})^{43/8}}{n^{3/8}}\exp\left(-\frac{3\omega_n}{8}\right) \text{ for Case I.}
\end{equation}
Here~\(D_0 > 0\) is as in~(\ref{sj_sparse}).
Similarly for case~\((II),\) there are at least~\(5\) squares star adjacent to~\(A.\) Choosing exactly~\(5\) such squares and using~(\ref{sj_sparse}) with~\(q = 5,\) we have
\begin{equation}\label{case2}
\mathbb{P}(J(A)) \leq D_0 \frac{(\log{n})^{29/8}}{n^{5/8}}\exp\left(-\frac{5\omega_n}{8}\right) \text{ for Case II.}
\end{equation}
Finally, for case~\((III),\) there are~\(8\) squares star adjacent to~\(A\) and so using~(\ref{sj_sparse}) with~\(q = 8,\) we have
\begin{equation}\label{case3}
\mathbb{P}(J(A)) \leq D_0\frac{\log{n}}{n}e^{-\omega_n} \text{ for Case III.}
\end{equation}
As before, let~\(N_j\) be the number of squares in~\(\{S_k\}\) satisfying Case~\((j)\) for~\(j \in \{I,II,III\}.\) Using the estimates for~\(N_I,N_{II}\) and~\(N_{III}\) in~(\ref{nj_est}), we have from~(\ref{case1}),~(\ref{case2}) and~(\ref{case3}) that
\begin{eqnarray}
\mathbb{P}\left(J_n \right) &\leq& 4D_0 \frac{(\log{n})^{43/8}}{n^{3/8}}\exp\left(-\frac{5\omega_n}{8}\right) + \sqrt{n}D_0 \frac{(\log{n})^{29/8}}{n^{5/8}}\exp\left(-\frac{3\omega_n}{8}\right) \nonumber\\
&&\;\;\;\;\;\;\;\;\;+\;\;\;\frac{8n}{\log{n}} D_0\frac{\log{n}}{n}e^{-\omega_n} \nonumber\\
&\leq& D e^{-\omega_n}\label{fi_est}
\end{eqnarray}
for all~\(n\) large and some constant~\(D > 0.\) The estimate~(\ref{fi_est}) is obtained using estimates~(\ref{ei_est1}) and~(\ref{ei_est2}). This proves~(\ref{fi_main_est}).~\(\qed\)

\subsection*{Constructing the Hamiltonian cycle}
Define the event
\begin{equation}\label{qn_def}
H_n = F_n \bigcap I_n^c \bigcap J_n^c
\end{equation}
where~\(F_n\) is as defined in~(\ref{q_n_def2}), the events~\(I_n\) and~\(J_n\) are as in~(\ref{ei_def}) and~(\ref{fi_def}), respectively. From~(\ref{q_n_est2}),~(\ref{ei_main_est}) and~(\ref{fi_main_est}), we have that
\begin{equation}\label{hn_est}
\mathbb{P}(H_n) \geq 1-\frac{1}{n^{9}} - 2De^{-\omega_n} \geq 1-3De^{-\omega_n}
\end{equation}
for all~\(n\) large. The final estimate is true since~\(\frac{\omega_n}{\log{\log{n}}} \longrightarrow 0\) as~\(n \rightarrow \infty.\) If the event~\(H_n\) occurs, then there is a backbone~\({\cal B}\) containing dense squares. The backbone~\({\cal B}\) is a dense star connected component and since the event~\(I_n^c\) occurs, there is no other dense star connected component. Also since~\(J_n^c\) also occurs, every sparse square is star adjacent to some dense square in~\({\cal B}.\)

We obtain the desired Hamiltonian cycle as in the case of long cycles~in Section~\ref{pf_long_cycle_rgg}. Let~\({\cal B} = \{W_i\}_{1 \leq i \leq t}\) be the set of dense squares in the backbone~\({\cal B}\) and for~\(1 \leq i \leq t,\) let~\(\eta_i\) be a cycle of edges containing all vertices in the square~\(W_i.\) As in the proof of~(\ref{eq_lg1_rgg}) in Theorem~\ref{long_rgg}, we obtain the cycle~\(\tau({\cal B})\) containing all vertices of~\({\cal B}.\) 

We now iteratively expand the cycle~\(\chi_0 := \tau({\cal B})\) by considering sparse squares attached to dense squares in~\({\cal B}.\) More precisely, let~\(\{Z_1,\ldots,Z_b\} \subset \{S_j\}\) be the set of all sparse squares. For~\(1 \leq j \leq b,\) let~\(\xi_j, 1 \leq j \leq b\) be any path containing all vertices in the square~\(Z_j.\) We iteratively construct a sequence of cycles~\(\{\chi_i\}_{1 \leq i \leq b}\) using the paths~\(\{\xi_j\}_{1 \leq j \leq b}.\)

Fix~\(1 \leq j \leq b.\) The cycle~\(\chi_{j-1}\) satisfies the following properties.\\
\((c1)\) The cycle~\(\chi_{j-1}\) contains all edges from~\(\xi_i, 1 \leq i \leq j-1\) and all edges from~\(\tau({\cal B})\) not removed so far in the iteration process.\\
\((c2)\) The square~\(Z_j\) is star adjacent to some dense square~\(W_j \in {\cal B}.\) There is at least one edge~\(g_j \in \chi_{j-1} \cap {\cal T}({\cal B})\) having both endvertices in~\(W_j.\)\\
\emph{Proof of~\((c1)-(c2)\) for~\(i=0\)}: Since~\(Z_j \notin {\cal B}\) is star adjacent to~\(W_j,\) at most~\(7\) squares of~\({\cal B}\) are star adjacent to~\(W_j.\) But the square~\(W_j\) contains at least~\(8\) vertices and so the cycle~\(\eta_j\) contained in~\(W_j\) contains at least~\(8\) edges. In each iteration in the proof of~(\ref{eq_lg1_rgg}) in Theorem~\ref{long_rgg}, at most one edge from each cycle~\(\eta_j\) was removed and so there is at least one edge of~\(\eta_j\) belonging to~\(\tau({\cal B}).\)~\(\qed\)

We remove the edge~\(g_{j}\) from the cycle~\(\chi_{j-1}\) to get a path~\(\chi_{j-1}\setminus \{g_j\}\) with endvertices~\(a_j\) and~\(b_j\) contained in the square~\(W_j.\) Similarly, let~\(c_j\) and~\(d_j\) denote the endvertices of the path~\(\xi_j\) contained in the square~\(Z_j.\) Since~\(Z_j\) is star adjacent to~\(W_j,\) the vertices~\(a_j\) and~\(c_j\) are joined by an edge. Similarly~\(b_j\) and~\(d_j\) are joined by an edge. We then merge the two paths to get the cycle~\[\chi_{j} = \left(\chi_{j-1} \setminus \{g_j\}\right) \bigcup \xi_j.\] As above, the new cycle~\(\chi_j\) also satisfies properties~\((c1)-(c2).\) We then repeat the above procedure until all the paths~\(\{\xi_j\}_{1 \leq j \leq b}\) have been merged. The above procedure continues for~\(b\) steps and the final cycle~\(\chi_b\) is the desired Hamiltonian cycle.~\(\qed\)

\subsection*{Acknowledgement}
I thank Professors Rahul Roy and Federico Camia for crucial comments and for my fellowships.

\bibliographystyle{plain}

\end{document}